\documentclass[11pt,a4paper]{article}

\usepackage{iftex}
\ifPDFTeX
  \usepackage[utf8]{inputenc}
  \usepackage[T1]{fontenc}
  \usepackage{lmodern}
\else
  \usepackage{fontspec}
\fi
\usepackage[english]{babel}
\usepackage[leqno]{amsmath}
\usepackage{amssymb}
\usepackage{mathrsfs}
\usepackage[normalem]{ulem}
\usepackage[a4paper,margin=2.6cm,marginparwidth=1.6cm,marginparsep=4mm]{geometry}

\newcommand{\EE}{\mathbb{E}}
\newcommand{\PP}{\mathbb{P}}
\newcommand{\RR}{\mathbb{R}}
\newcommand{\DD}{\mathbb{D}}
\newcommand{\Etrib}{\mathscr{E}}
\newcommand{\Ftrib}{\mathscr{F}}
\newcommand{\Mcal}{\mathscr{M}}
\newcommand{\Dsc}{\mathscr{D}}
\newcommand{\lsup}[1]{{}^{(#1)}\!}
\newcommand{\cro}[1]{\langle #1\rangle}
\newcommand{\pgo}[1]{\marginpar{\footnotesize\itshape p.~#1}}
\newenvironment{enonce}[1]%
  {\par\medskip\noindent\uline{#1}.\kern.6em\itshape}%
  {\par\medskip}
\newcommand{\dem}{\par\medskip\noindent\uline{Proof}.\kern.6em}
\newcommand{\findem}{\hfill$\square$\par\medskip}
\newcommand{\rem}[1]{\par\medskip\noindent\uline{#1}.\kern.6em}

\begin{document}

\noindent{\small\itshape English translation of N.~Bouleau,
``Décomposition de l'énergie par niveau de potentiel,''  (Colloque de théorie du potentiel Jacques Deny, Orsay, June 1983, Lecture Notes in Mathematics 1096, Springer,
1984, pp.~149--172.) }

\vspace{1em}\hrule\vspace{2.5em}

\pgo{149}%
\noindent
\begin{tabular}{@{}l@{}}
Colloquium on Potential Theory---\\
in honor of Jacques Deny\\
\hspace*{1.2em}--- Orsay 1983 ---
\end{tabular}

\vspace{3em}

\begin{center}
\uline{ENERGY DECOMPOSITION BY POTENTIAL LEVEL}\\[1.2em]
{\itshape Nicolas BOULEAU}
\end{center}

\bigskip

We study various extensions of the following result:

\begin{enonce}{PROPOSITION 1}
Let $m$ be Lebesgue measure on $\RR^d$, and let $u\in
L^2_{\mathrm{loc}}(\RR^d)$ have partial derivatives in distributions sense in
$L^2(\RR^d)$. Then the image under $u$ of the measure
$\operatorname{grad}^2u\, .\, m$ is absolutely continuous with respect
to Lebesgue measure on $\RR$.
\end{enonce}

In the first part we study excessive functions for a Markov process
that are continuous semimartingales along the sample paths. The
property then follows from the theory of local times of semimartingales.
We next treat the case of Dirichlet spaces and show that, if $\tilde u$
is a quasi-continuous version of a function $u$ belonging to a regular
Dirichlet space, then the image under $\tilde u$ of the local
energy measure of $u$ is absolutely continuous with respect to Lebesgue
measure. Finally, in a third part, we establish the occupation-time
density property for certain Dirichlet processes.

\pgo{150}%
\bigskip
\noindent
I. \uline{EXCESSIVE FUNCTIONS OF MARKOV PROCESSES}.

\medskip
a) Let $(\Omega,X,\PP_x)$ be a right process with state space $E$ and
canonical $\sigma$-fields $(\Ftrib_t)_{t\ge 0}$ (cf.~[10]).

Let $u$ be a finite excessive function that is the potential of a
continuous adapted increasing additive functional $A$,
\begin{equation*}\tag{1}
u(x) = \EE_x\,A_\infty
\end{equation*}
and that can be written along the paths as
\begin{equation*}\tag{2}
u(X_t) = u(X_0) + M_t - A_t,
\end{equation*}
where the martingale starting from zero,
\[
M_t = \EE[\,A_\infty\mid\Ftrib_t\,] - \EE[\,A_\infty\mid\Ftrib_0\,],
\]
is assumed continuous for simplicity.

To $u$ one may associate the \uline{energy process} $EP(u)$ defined by
\begin{equation*}\tag{3}
EP(u)_t = \tfrac12\bigl[\,u^2(X_0) + \cro{M,M}_t\,\bigr],
\end{equation*}
the \uline{energy function} $EF(u)$ given by
\begin{equation*}\tag{4}
EF(u)(x) = \tfrac12\bigl[\,u^2(x) + \EE_x\cro{M,M}_\infty\,\bigr],
\end{equation*}
and, after specifying a positive measure $\theta$ on $E$, the
\uline{energy number} $EN_\theta(u)$:
\begin{equation*}\tag{5}
EN_\theta(u) = \tfrac12\bigl[\,\langle\theta,u^2\rangle
  + \EE_\theta\cro{M,M}_\infty\,\bigr].
\end{equation*}

We then have
\[
EF(u)(x) = \tfrac12\,\EE_x\bigl[\,u(X_0)+M_\infty\,\bigr]^2
         = \tfrac12\,\EE_x A_\infty^2,
\]
which, by the energy formula (cf.~[6]), gives
\[
EF(u)(x) = \EE_x\int_0^\infty u(X_s)\,dA_s.
\]
Thus, if $S^u$ denotes the excessive kernel associated with $u$, defined
by
\[
S^u f = \EE_.\int_0^\infty f(X_s)\,dA_s,
\]
then
\begin{equation*}\tag{6}
\underline{EF(u) = S^u u}.
\end{equation*}

\pgo{151}%
\medskip
b) \uline{Representation of the energy as an integral of potentials
supported on the level sets of $u$}.

Write $Y$ for the semimartingale $u(X_t)=u(X_0)+M_t-A_t$, and let
$L^a_t$, $a\in\RR$, be the local time of $Y$ at $a$. It is a continuous
increasing process satisfying the Meyer--Tanaka formula (cf.~[13]):
\begin{equation*}\tag{7}
(u\wedge a)(X_t) = (u\wedge a)(X_0)
 + \int_0^t \mathbf{1}_{\{u(X_s)<a\}}(dM_s-dA_s) - \tfrac12 L^a_t,
\end{equation*}
and the occupation-time density property holds:
\begin{equation*}\tag{8}
\int_0^t g(Y_s)\,d\cro{M,M}_s = \int_\RR g(a)\,L^a_t\,da
 \qquad \PP_x\text{-a.s. for every }x,
\end{equation*}
for every positive Borel function $g$.

It follows from (7) that the excessive kernel $S^{u\wedge a}$ associated
with the excessive function $u\wedge a$ is given by
\[
S^{u\wedge a}(f)(x) = \EE_x\Bigl[\,\int_0^\infty
 f(X_s)\,\mathbf{1}_{\{u(X_s)<a\}}\,dA_s
 + \tfrac12\int_0^\infty f(X_s)\,dL^a_s\,\Bigr].
\]
Since the local time $L^a_t$ is carried by
$\{(\omega,t):u(X_t(\omega))=a\}$, taking
$f=\mathbf{1}_{\{u=a\}}$ gives
\begin{equation*}\tag{9}
S^{u\wedge a}\,\mathbf{1}_{\{u=a\}} = \tfrac12\,\EE_x\,L^a_\infty,
\end{equation*}
and therefore, by (8) and (4),
\begin{equation*}\tag{10}
\underline{EF(u) = S^u u = \int_0^\infty
 \bigl(S^{u\wedge a}\,\mathbf{1}_{\{u=a\}}\bigr)\,da + \tfrac12 u^2}.
\end{equation*}

Assume that there exists a reference measure $\xi$, chosen excessive
and $\sigma$-finite. With every continuous additive functional $B_t$
one may then associate (cf.~[15]) a positive $\sigma$-finite measure
$\mu_B$, which does not charge semipolar sets, given by
\begin{equation*}\tag{11}
\mu_B(f) = \lim_{t\downarrow 0}\frac1t\,
 \EE_\xi\int_0^t f(X_s)\,dB_s,
 \qquad f\text{ positive and measurable}.
\end{equation*}

The measure $\nu_u$ associated by (11) with $\cro{M,M}_t$ will be called
the \uline{energy measure} of $u$. Similarly, to the additive functionals
$L^a_t$ correspond measures $\mu^a_u$ supported on $\{u=a\}$; these
provide a disintegration of $\nu_u$ under the mapping $x\mapsto u(x)$.

\pgo{152}%
Indeed, by (8),
\[
\int f\circ u(x)\,d\nu_u(x)
 = \int_\RR\Bigl(\,\int f(x)\,d\mu^a_u(x)\Bigr)\,da,
 \qquad f\text{ positive and measurable},
\]
and the image under $u$ of the restriction of $\nu_u$ to any set
on which it is finite is absolutely continuous with respect to Lebesgue
measure.

\medskip
c) Under the same assumptions as in a), let $\varphi$ be an increasing
concave function from $\RR_+$ into $\RR_+$ that vanishes at zero. Then
$\varphi\circ u$ satisfies the same assumptions as $u$. By the
change-of-variables formula for local times (cf.~[3]), if
$b=\varphi(a)$, then
\[
L^b(\varphi\circ Y) = \varphi_g'(a)\,L^a(Y)
 \qquad \forall a>0,
\]
where $\varphi'_g$ is the left derivative of $\varphi$. It follows from
(9) and (10) that
\begin{align*}
EF(\varphi\circ u)
 &= \tfrac12(\varphi\circ u)^2
   + \tfrac12\int_0^\infty \EE_.\bigl[\,L^b_\infty(\varphi\circ Y)\,\bigr] db\\
 &= \tfrac12(\varphi\circ u)^2
   + \tfrac12\int_0^\infty \varphi_g'^{\,2}(a)\,
     \EE_.\bigl[\,L^a_\infty\,\bigr] da.
\end{align*}
Hence the following change-of-variables formula for the energy:
\begin{equation*}\tag{12}
\underline{EF(\varphi\circ u) = \tfrac12(\varphi\circ u)^2
 + \int_0^\infty \varphi'^{\,2}(a)\,
   \bigl(S^{u\wedge a}\,\mathbf{1}_{\{u=a\}}\bigr)\,da},
\end{equation*}
where any version of the Lebesgue derivative of $\varphi$ may be used.
This formula could be extended under less restrictive assumptions, but
it is expressed more naturally in the setting of Dirichlet spaces, as
we shall now see.

\bigskip
\noindent
II. \uline{THE CASE OF DIRICHLET SPACES}

\medskip
We follow Fukushima's presentation [8], with the following notation and
assumptions.

Let $E$ be a locally compact space with a countable base, equipped with
its Borel $\sigma$-field $\Etrib$, and let $m$ be a positive
\pgo{153}%
$\sigma$-finite measure with full support $E$. Let $\Phi$ be a positive
symmetric bilinear form defined on the vector subspace $\Dsc\Phi$ of
$L^2(E,\Etrib,m)$. We assume that:

\smallskip
\begin{itemize}
\item[$\bullet$] $\Dsc\Phi\cap C_K(E)$ is dense in $\Dsc\Phi$ for the
norm associated with the form $\Phi_1(u,v)=\Phi(u,v)+(u,v)$, where
$(\cdot,\cdot)$ is the inner product of $L^2(m)$, and is dense in
$C_K(E)$ for the natural topology of $C_K(E)$;
\item[$\bullet$] $\Phi$ is closed, i.e. $\Dsc\Phi$ is complete for the
metric associated with $\Phi_1$;
\item[$\bullet$] the unit contraction operates, i.e. $u\in\Dsc\Phi$ and
$v=(u\vee0)\wedge1$ imply $v\in\Dsc\Phi$ and
$\Phi(v,v)\le\Phi(u,u)$.
\end{itemize}
\smallskip

In other words, $\Phi$ is a regular Dirichlet form on $L^2(m)$; we
denote by $\DD=\Dsc\Phi$ the associated Dirichlet space.

There then exists a Hunt process $(\Omega,\Ftrib_t,X_t,P_x)$ with state
space $E\cup\{\delta\}$ whose transition function $P_t$ is
$m$-symmetric,
\[
(P_tu,v)=(u,P_tv),
\]
for all positive measurable $u,v$ on $E$, and which defines a strongly
continuous semigroup on $L^2(m)$ such that
\begin{gather*}
\DD = \Bigl\{u\in L^2(m):\lim_{t\downarrow0}\frac1t\,(u-P_tu,u)
 <+\infty\Bigr\},\\
\Phi(u,u)=\lim_{t\downarrow0}\uparrow\frac1t\,(u-P_tu,u)
 \qquad \forall u\in\DD.
\end{gather*}

If $u\in\DD$, $\tilde u$ denotes a quasi-continuous version of $u$.
``Quasi-everywhere'' means outside an $m$-polar set.

There exists ([8], Lemma 4.5.2) a positive $\sigma$-finite measure $k$
on $E$ that does not charge $m$-polar sets and such that
\begin{equation*}\tag{13}
\forall u\in\DD \quad \lim_{t\downarrow0}\frac1t\,
\EE_m\bigl[\,u^2(X_0)\bigl(\mathbf{1}_E(X_0)-\mathbf{1}_E(X_t)\bigr)\bigr]
 = \langle\tilde u^2,k\rangle,
\end{equation*}
and, for every $u\in\DD$,
\begin{equation*}\tag{14}
\Phi(u,u) = \lim_{t\downarrow0}\frac1{2t}\,
\EE_m\bigl[\,(u(X_t)-u(X_0))^2\,\bigr]
 + \lim_{t\downarrow0}\frac1t\,
\EE_m\bigl[\,u^2(X_0)\bigl(\mathbf{1}_E(X_0)-\mathbf{1}_E(X_t)\bigr)\bigr].
\end{equation*}

Recall the following result of Fukushima.

\pgo{154}%
\begin{enonce}{PROPOSITION 2}
Let $u\in\DD$. Then there exists an $m$-polar set $N(u)$ such that, for
every $x$ outside $N(u)$,
\begin{equation*}\tag{15}
\tilde u(X_t) = \tilde u(X_0) + \lsup{u}M_t + \lsup{u}A_t
\end{equation*}
up to a $\PP_x$-evanescent set, where $\lsup{u}M_t$ is a martingale
additive functional satisfying
$\EE_x[\,\lsup{u}M_t^2\,]<\infty$ for every $t$, and where
$\lsup{u}A_t$ is an additive functional of zero energy:
\[
\lim_{t\downarrow0}\frac1{2t}\,\EE_m[\,\lsup{u}A_t^2\,]=0,
\]
and
\begin{equation*}\tag{16}
\lim_{t\downarrow0}\frac1{2t}\,
\EE_m\bigl[\,(u(X_t)-u(X_0))^2\,\bigr]
 = \lim_{t\downarrow0}\frac1{2t}\,\EE_m[\,\lsup{u}M_t^2\,]
 = \sup_{t>0}\frac1{2t}\,\EE[\,\lsup{u}M_t^2\,].
\end{equation*}
\end{enonce}

Let $\Mcal$ be the set of processes $M$ which, for quasi-every $x$, are
under $\PP_x$ square-integrable martingale additive functionals in the
wide sense and satisfy
\[
e(M)=\sup_{t>0}\frac1{2t}\,\EE_m[\,M_t^2\,]<+\infty.
\]

Equipped with the inner product associated with $e$, $\Mcal$ is a
Hilbert space. By Doob's inequality, the subspace of $\Mcal$ consisting
of continuous martingales is closed in $\Mcal$ ([8], Theorem 5.2.1).
We write $M^c$ for the projection of $M$ onto this subspace and
$M^d=M\ominus M^c$.

The decomposition of $u\in\DD$,
\begin{equation*}\tag{17}
\tilde u(X_t) = \tilde u(X_0) + \lsup{u}M^c_t + \lsup{u}M^d_t
 + \lsup{u}A_t,
\end{equation*}
valid up to a $\PP_x$-evanescent set for quasi-every $x$, will be called
the \uline{canonical decomposition} of $u$.

\begin{enonce}{LEMMA 3}
Let $u\in\DD$, let $f$ be a compactly supported $C^1$ function, and set
$F(x)=\int_0^x f(y)\,dy$.

Then $F\circ u\in\DD$, and the canonical decomposition of $F\circ u$
is
\begin{equation*}\tag{18}
F\circ\tilde u(X_t) = F\circ\tilde u(X_0)
 + \int_0^t f\circ\tilde u(X_s)\,d\lsup{u}M^c_s
 + \lsup{F\circ u}M^d_t + \lsup{F\circ u}A_t
\end{equation*}
under $P_x$ for quasi-every $x$.
\end{enonce}

\pgo{155}%
\dem
Set $u_n=nU_{n+1}u=U_1v_n$, with $v_n=n(u-nU_{n+1}u)$, where
$(U_p)_{p>0}$ is the resolvent of $(P_t)$.

The canonical decomposition of $u_n$ is
\[
u_n(X_t) = u_n(X_0) + \lsup{u_n}M^c_t + \lsup{u_n}M^d_t
 + \int_0^t\bigl(u_n(X_s)-v_n(X_s)\bigr)\,ds,
\]
valid under $\PP_x$ for every $x$.

Hence, by Itô's formula,
\begin{align*}
F\circ u_n(X_t) ={}& F\circ u_n(X_0)
 + \int_0^t f\circ u_n(X_s)\,d\lsup{u_n}M^c_s
 + \int_0^t f\circ u_n(X_{s-})\,d\lsup{u_n}M^d_s\\
&+ \int_0^t f\circ u_n(X_s)\bigl(u_n(X_s)-v_n(X_s)\bigr)\,ds
 + \frac12\int_0^t f'\circ u_n(X_s)\,
   d\cro{\lsup{u_n}M^c,\lsup{u_n}M^c}_s\\
&+ \sum_{0<s\le t}\bigl[\,F\circ u_n(X_s)-F\circ u_n(X_{s-})
 - f\circ u_n(X_{s-})\bigl(u_n(X_s)-u_n(X_{s-})\bigr)\bigr].
\end{align*}

The canonical decomposition of $F\circ u_n$ is therefore of the form
\begin{equation*}\tag{19}
F\circ u_n(X_t) = F\circ u_n(X_0)
 + \int_0^t f\circ u_n(X_s)\,d\lsup{u_n}M^c_s
 + \lsup{F\circ u_n}M^d_t + \lsup{F\circ u_n}A_t.
\end{equation*}

As $n\to\infty$, $u_n\to u$ in $(\DD,\Phi_1)$, and hence
$F\circ u_n\to F\circ u$ in $(\DD,\Phi_1)$, because
$F/\|f\|_\infty$ is a normal contraction and is therefore continuous
on $(\DD,\Phi_1)$ (cf.~[1]). It follows ([8], Theorem 5.2.2) that
\[
\lsup{u_n}M_t \to \lsup{u}M_t
 \quad\text{and}\quad
\lsup{F\circ u_n}M_t \to \lsup{F\circ u}M_t
 \qquad\text{in }(\Mcal,e),
\]
and hence also that the continuous and purely discontinuous parts
satisfy
\begin{equation*}\tag{20}
\lsup{u_n}M^c_t \to \lsup{u}M^c_t,\quad
\lsup{F\circ u_n}M^c_t \to \lsup{F\circ u}M^c_t,\quad
\lsup{F\circ u_n}M^d_t \to \lsup{F\circ u}M^d_t
\end{equation*}
in $(\Mcal,e)$.

Choose a subsequence such that $u_m\to\tilde u$ quasi-everywhere. From
(20), and from the fact that $f$ is bounded and continuous, we obtain
\[
\int_0^t f\circ u_m(X_s)\,d\lsup{u_m}M^c_s
 \to \int_0^t f\circ\tilde u(X_s)\,d\lsup{u}M^c_s
\]
in $(\Mcal,e)$.
\pgo{156}%
Thus
\[
\lsup{F\circ u}M^c_t = \int_0^t f\circ\tilde u(X_s)\,d\lsup{u}M^c_s,
 \qquad\text{which is (18).}
\]
\findem

Since $X_t$ is a Hunt process, there exists a canonical continuous
increasing additive functional $\xi_s$. Let $(N(x,dy),\xi_s)$ be the
associated Lévy system. Recall that, if $h(x,y)$ is a positive
measurable function on $E\times E$ that vanishes on the diagonal, then
the predictable projection of the random measure
\[
\sum_{s>0}h(X_s,X_{s-})\,\varepsilon_s
\]
is the measure
\[
\Bigl[\,\int h(y,X_s)\,N(X_s,dy)\Bigr]d\xi_s.
\]

We then have:

\begin{enonce}{COROLLARY 4}
Let $u\in\DD$, let $f$ be a compactly supported $C^1$ function, and set
$F(x)=\int_0^x f(y)\,dy$. Then
\begin{equation*}\tag{21}
\begin{split}
\cro{\lsup{F\circ u}M,\lsup{F\circ u}M}_t
 ={}& \int_0^t f^2\circ\tilde u(X_s)\,
     d\cro{\lsup{u}M^c,\lsup{u}M^c}_s\\
 &+ \int_0^t\!\!\int
   \bigl(F\circ\tilde u(y)-F\circ\tilde u(X_s)\bigr)^2
   N(X_s,dy)\,d\xi_s
\end{split}
\end{equation*}
under $\PP_x$ for quasi-every $x$.
\end{enonce}

\dem
It follows from (18) that
\[
\Delta\lsup{F\circ u}M^d_t
 = F\circ\tilde u(X_t)-F\circ\tilde u(X_{t-}).
\]
The process $\cro{\lsup{F\circ u}M^d,\lsup{F\circ u}M^d}_t$, the dual
predictable projection of
$\sum_{0<s\le t}(\Delta\lsup{F\circ u}M^d_s)^2$, is therefore equal to
\[
\int_0^t\!\!\int
 \bigl(F\circ\tilde u(y)-F\circ\tilde u(X_s)\bigr)^2N(X_s,dy)\,d\xi_s,
\]
which proves the corollary.
\findem

\pgo{157}%
Let $u\in\DD$, and denote by $\mu_u$ the \uline{local energy measure}
of $u$, defined by
\[
\langle\mu_u,h\rangle = \lim_{t\downarrow0}\frac1{2t}\,
 \EE_m\int_0^t h(X_s)\,d\cro{\lsup{u}M^c,\lsup{u}M^c}_s
\]
for every positive $\Etrib$-measurable function $h$.

The measure $\mu_u$ does not charge $m$-polar sets, and
\[
\|\mu_u\| \le \lim_{t\downarrow0}\frac1{2t}\,
 \EE_m[\,\lsup{u}M_t^2\,] \le \Phi(u,u)<+\infty.
\]
We may therefore also define $\nu_u$ as the image of $\mu_u$ under
$\tilde u$. By (18),
\begin{equation*}\tag{22}
\lim_{t\downarrow0}\frac1{2t}\,
 \EE_m\bigl[\,\bigl(\lsup{F\circ u}M^c_t\bigr)^2\,\bigr]
 = \langle\mu_u,f^2\circ\tilde u\rangle
 = \int f^2(y)\,d\nu_u(y).
\end{equation*}

Now let $g$ be a bounded Borel function and set
$G(x)=\int_0^x g(y)\,dy$. Let $g_n$ be compactly supported $C^1$
functions such that $|g_n|\le\|g\|_\infty$ and
$g_n\to g$ in $L^2\bigl(\nu_u+\frac{dx}{1+x^2}\bigr)$.

a) First, $g_n\to g$ in $L^1_{\mathrm{loc}}$, and therefore
$G_n=\int_0^{\,\cdot}g_n(y)\,dy\to G$ everywhere. Since
$|G_n(y)|\le\|g\|_\infty|y|$, it follows that
\begin{equation*}\tag{23}
G_n\circ u\to G\circ u \qquad\text{in }L^2(m).
\end{equation*}

b) We show that $G_n\circ u$ is a Cauchy sequence for $\Phi$. By
(13), (14), and (16),
\begin{align*}
\Phi(G_p\circ u-G_q\circ u\,,\,G_p\circ u-G_q\circ u)
 ={}& \lim_{t\downarrow0}\frac1{2t}\,
   \EE_m\bigl[\,\bigl(\lsup{G_p\circ u-G_q\circ u}M_t\bigr)^2\,\bigr]\\
 &+ \bigl\langle(G_p\circ\tilde u-G_q\circ\tilde u)^2,k\bigr\rangle.
\end{align*}

$\bullet$ Since
$|G_p\circ\tilde u-G_q\circ\tilde u|
 \le2\|g\|_\infty\tilde u$, and
$\langle\tilde u^2,k\rangle\le\Phi(u,u)<+\infty$, by
$G_p-G_q\to0$, we have
\[
\langle(G_p\circ\tilde u-G_q\circ\tilde u)^2,k\rangle\to0
 \qquad\text{as }p,q\uparrow\infty.
\]

$\bullet$ Let $\eta$ be the positive $\sigma$-finite measure, not
charging $m$-polar sets, associated with the canonical additive
functional $\xi_t$. Thus
\[
\langle\eta,h\rangle = \lim_{t\downarrow0}\frac1t\,
 \EE_m\int_0^t h(X_s)\,d\xi_s,
 \qquad h\text{ positive and measurable}.
\]

\pgo{158}%
Corollary 4 gives
\begin{align*}
\lim_{t\downarrow0}\frac1{2t}\,
 \EE_m\bigl[\,\lsup{G_p\circ u-G_q\circ u}M_t^2\,\bigr]
 ={}& \langle(g_p-g_q)^2\circ\tilde u,\mu_u\rangle\\
 &+ \tfrac12\bigl\langle\textstyle\int
   \bigl[(G_p-G_q)\circ\tilde u(y)
   -(G_p-G_q)\circ\tilde u(x)\bigr]^2N(x,dy),\eta(dx)\bigr\rangle.
\end{align*}

The first term, which equals $\langle(g_p-g_q)^2,\nu_u\rangle$, tends
to zero because $g_p\to g$ in $L^2(\nu_u)$.

For the second term, we have
\[
\bigl|(G_p-G_q)\circ\tilde u(y)-(G_p-G_q)\circ\tilde u(x)\bigr|
 \le2\|g\|_\infty|\tilde u(y)-\tilde u(x)|.
\]
Moreover,
\begin{align*}
\int\!\!\int(\tilde u(y)-\tilde u(x))^2N(x,dy)\,\eta(dx)
 &= \lim_{t\downarrow0}\frac1t\,\EE_m\int_0^t\!\!\int
   (\tilde u(y)-\tilde u(X_s))^2N(X_s,dy)\,d\xi_s\\
 &= \lim_{t\downarrow0}\frac1t\,\EE_m\Bigl[\sum_{0<s\le t}
   \bigl(u(X_s)-u(X_{s-})\bigr)^2\Bigr]\\
 &= \lim_{t\downarrow0}\frac1t\,\EE_m\Bigl[\sum_{0<s\le t}
   \bigl(\Delta\lsup{u}M_s\bigr)^2\Bigr]\\
 &\le \lim_{t\downarrow0}\frac1t\,\EE_m[\,\lsup{u}M_t^2\,]
 \le \Phi(u,u)<+\infty.
\end{align*}
The dominated convergence theorem therefore applies, proving that
$G_n\circ u$ is a Cauchy sequence for $\Phi$.

c) By a) and b), since $\Phi$ is closed,
$G_n\circ u\to G\circ u$ in $(\DD,\Phi_1)$.

As in the proof of Lemma 3, this implies that
\[
\lsup{G_n\circ u}M^c_t\to\lsup{G\circ u}M^c_t
 \qquad\text{in }(\Mcal,e).
\]
Since $g_n\to g$ in $L^2(\nu_u)$ also implies
\[
\int_0^t g_n\circ\tilde u(X_s)\,d\lsup{u}M^c_s
 \to \int_0^t g\circ\tilde u(X_s)\,d\lsup{u}M^c_s
 \qquad\text{in }(\Mcal,e),
\]
it follows that the canonical decomposition of $G\circ u$ is
\[
G\circ\tilde u(X_t)=G\circ\tilde u(X_0)
 +\int_0^t g\circ\tilde u(X_s)\,d\lsup{u}M^c_s
 +\lsup{G\circ u}M^d_t+\lsup{G\circ u}A_t.
\]
Consequently,
\pgo{159}%
\begin{align*}
\Phi(G\circ u,G\circ u)
 &= \lim_{t\downarrow0}\frac1{2t}\,
   \EE_m\bigl[\,\lsup{G\circ u}M_t^2\,\bigr]
   +\langle G^2\circ\tilde u,k\rangle\\
 &= \langle g^2,\nu_u\rangle
   +\int\!\!\int\bigl(G\circ\tilde u(y)-G\circ\tilde u(x)\bigr)^2
     N(x,dy)\,\eta(dx)
   +\langle G^2\circ\tilde u,k\rangle.
\end{align*}

We have therefore proved:

\begin{enonce}{THEOREM 5}
Let $g$ be a bounded Borel function and let
$G(x)=\int_0^x g(y)\,dy$. Then, for every $u\in\DD$,
$G\circ u\in\DD$, and
\begin{equation*}\tag{24}
\Phi(G\circ u,G\circ u)=\langle g^2,\nu_u\rangle
 +\int\!\!\int\bigl(G\circ\tilde u(y)-G\circ\tilde u(x)\bigr)^2
   N(x,dy)\,\eta(dx)
 +\langle G^2\circ\tilde u,k\rangle.
\end{equation*}
The canonical decomposition of $G\circ u$ is
\begin{equation*}\tag{25}
G\circ\tilde u(X_t)=G\circ\tilde u(X_0)
 +\int_0^t g\circ\tilde u(X_s)\,d\lsup{u}M^c_s
 +\lsup{G\circ u}M^d_t+\lsup{G\circ u}A_t
\end{equation*}
under $\PP_x$ for quasi-every $x$.
\end{enonce}

For a function of one variable, this extends a formula of Le Jan [11].

\begin{enonce}{COROLLARY 6}
The measure $\nu_u$, the image under $\tilde u\in\DD$ of the
local energy measure $\mu_u$, is absolutely continuous with respect to
Lebesgue measure.\footnote{This result is a precursor to the EID (energy image density) property of Bouleau and Hirsch, demonstrated on Wiener space in 1985 (Bouleau, N. and Hirsch, F. Propriétés d’absolue continuité dans les espaces de Dirichlet et application aux équations différentielles stochastiques, Sém. de Probabilités XX, 1984/85, pp. 131–161, LNM 1204, Springer 1986, and in the book Bouleau, N. and Hirsch, F. Dirichlet forms and analysis on Wiener space, De Gruyter Stud. Math., 14,  1991). It was extended as a general conjecture in 1986 (Bouleau, N. and Hirsch, F. Formes de dirichlet générales et densité des variables aléatoires réelles sur l’espace de Wiener, J. Funct. Anal. 69 (1986), no. 2, pp. 229–259.).}
\end{enonce}

\dem
It is enough to take $g$ to be the indicator function of a
Lebesgue-null set in (24).

\rem{Remark 7}If the process $X_t$ is of Lebesgue type (i.e. if the
additive functional equal to $t$ is canonical), then, for every
$u\in\DD$,
\[
d\cro{\lsup{u}M^c,\lsup{u}M^c}_t\ll dt,
 \qquad\text{hence}\qquad \mu_u\ll m.
\]
Thus Corollary 6 holds with $u$ in place of $\tilde u$.

This is in particular the case for Brownian motion on $\RR^d$, which
proves Proposition 1 stated in the introduction. Proposition 1 can
also be proved directly by reducing to dimension one using Lemma 3.2 of
[7], and by showing that, if $u:\RR\to\RR$ is continuous and of bounded
variation, then the inverse image under $u$ of a Lebesgue-null set is
null for the measure $|du|$, $\;$ a property that is, incidentally,
nontrivial.

\pgo{160}%
\bigskip
\noindent
III. \uline{OCCUPATION-TIME DENSITY PROPERTIES}.

\medskip
A. We first remain in the setting of Dirichlet spaces.

\begin{enonce}{DEFINITION 8}\
We shall say that $u\in\DD$ satisfies quasi-everywhere the
occupation-time density property along the paths if, writing the
canonical decomposition of $u$ as
\[
\tilde u(X_t)=\tilde u(X_0)+\lsup{u}M^c_t+\lsup{u}M^d_t+\lsup{u}A_t
 \qquad \PP_x\text{-a.s. for every }x\notin N(u),
\]
the image of the measure
$d\cro{\lsup{u}M^c,\lsup{u}M^c}_s(\omega)$ on $[0,t]$ under the mapping
$s\mapsto\tilde u(X_s(\omega))$ is absolutely continuous with respect
to Lebesgue measure, $\PP_x$-a.s., for quasi-every $x$.
\end{enonce}

Suppose that $m$ is a reference measure, which is equivalent to saying
that the $m$-polar sets are polar (cf.~[8], Theorem 4.2.2). We are then
under the classical duality assumptions [2]. For $\alpha>0$, let
$u_\alpha(x,y)$ be the function symmetric and $\alpha$-excessive in
each variable such that
\[
U_\alpha f(x)=\int u_\alpha(x,y)f(y)\,dm(y)
 \qquad f\in\Etrib^+.
\]

A measure $\mu$ with bounded potential corresponds to an additive
functional $A_t$ such that
\[
U_\alpha(h.\mu)(x)=\int u_\alpha(x,y)h(y)\,d\mu(y)
 =\EE_x\int_0^\infty e^{-\alpha s}h(X_s)\,dA_s.
\]

It follows (cf.~[14], p.~765) that, if $\mu$ is a positive
$\sigma$-finite measure that does not charge polar sets, there
corresponds to it a positive homogeneous random measure $dA_t(\omega)$
such that
\begin{equation*}\tag{26}
U_\alpha(h.\mu)(x)=\EE_x\int_0^\infty e^{-\alpha s}h(X_s)\,dA_s
 \qquad \forall h\in\Etrib^+,
\end{equation*}
and therefore
\begin{equation*}\tag{27}
\langle h.\mu,U_\alpha g\rangle
 =\EE_{g.m}\int_0^\infty e^{-\alpha s}h(X_s)\,dA_s
 \qquad h,g\in\Etrib^+.
\end{equation*}

Let $u\in\DD$, and let $\tilde u$ be a Borel quasi-continuous version
of $u$. Disintegrating the local energy measure $\mu_u$ under the
mapping $x\mapsto\tilde u(x)$,
\pgo{161}%
we obtain measures $\mu^a_u$ such that
\begin{equation*}\tag{28}
\mu^a_u\text{ is supported on }(\tilde u=a),
 \qquad\text{and}\qquad \int\mu^a_u\,da=\mu_u.
\end{equation*}

Assume the following:

\medskip
\noindent
(29)\quad \uline{For Lebesgue-almost every} $a$, \uline{the measures}
$\mu^a_u$ \uline{do not charge polar sets}.

\medskip
This assumption does not depend on the version $\tilde u$, since
$\mu_u$ does not charge polar sets.

The measures $\mu^a_u$ then correspond to positive homogeneous random
measures $dA^a_t$ such that
\begin{equation*}\tag{30}
\langle h.\mu^a_u,U_\alpha g\rangle
 =\EE_{g.m}\int_0^\infty e^{-\alpha s}h(X_s)\,dA^a_s
 \qquad h,g\in\Etrib^+,\ \alpha>0.
\end{equation*}

The measure $\mu_u$ corresponds to a homogeneous random measure which,
by the definition of $\mu_u$, extends the measure
$d\cro{\lsup{u}M^c,\lsup{u}M^c}_s$, previously defined only
$\PP_x$-a.s. for quasi-every $x$. We continue to denote this measure by
$d\cro{\lsup{u}M^c,\lsup{u}M^c}_s$.

Thus, if $f$ is positive,
\begin{align*}
\EE_{g.m}\int_0^\infty e^{-\alpha s}h(X_s)
 f\circ\tilde u(X_s)\,d\cro{\lsup{u}M^c,\lsup{u}M^c}_s
 =\langle h\,.\,f\circ\tilde u\,.\,\mu_u,U_\alpha g\rangle,
\end{align*}
which, by (28) and then (30), equals
\begin{align*}
&=\int f(a)\,\langle h\,.\,\mu^a_u,U_\alpha g\rangle\,da\\
&=\EE_{g.m}\int_0^\infty\int_{a\in\RR}
 e^{-\alpha s}h(X_s)f(a)\,dA^a_s\,da.
\end{align*}

It follows that the $\alpha$-excessive functions
\[
\EE_.\int_0^\infty e^{-\alpha s}h(X_s)f\circ\tilde u(X_s)
 \,d\cro{\lsup{u}M^c,\lsup{u}M^c}_s
 \qquad\text{and}
\]
\[
\EE_.\int_0^\infty\int_{a\in\RR}
 e^{-\alpha s}h(X_s)f(a)\,dA^a_s\,da,
\]
which are equal $m$-almost everywhere, coincide everywhere.

Hence the random measures
\pgo{162}%
\[
f\circ\tilde u(X_s)\,d\cro{\lsup{u}M^c,\lsup{u}M^c}_s\,d\PP_x
 \qquad\text{and}\qquad
\int_{a\in\RR}f(a)\,dA^a_s\,da\,d\PP_x
\]
also coincide, for every $x$.

Now choose $x$ outside a polar set $N(u)$ so that
\[
\cro{\lsup{u}M^c,\lsup{u}M^c}_t<+\infty
 \qquad \PP_x\text{-a.s. for every }x\notin N(u),
\]
and let $f$ range over a countable dense subset of $C_K$. It follows
that $u$ satisfies quasi-everywhere the occupation-time density
property along the paths.

\rem{Remark 9}If $m$ is a reference measure and the empty set is the
only polar set, then condition (29) is trivially satisfied. This is in
particular the case for Brownian motion on $\RR$: for every
$u\in\DD=H^1(\RR)$,
\begin{gather*}
u(B_t)=u(B_0)+\int_0^t u'(B_s)\,dB_s+\lsup{u}A_t
 \qquad P_x\text{-a.s. for every }x,\quad\text{and}\\
d\cro{\lsup{u}M^c,\lsup{u}M^c}_s=u'^{\,2}(B_s)\,ds.
\end{gather*}

Then
\begin{equation*}\tag{31}
\varphi\longrightarrow\int_0^t\varphi\circ u(B_s)u'^{\,2}(B_s)\,ds
\end{equation*}
defines a measure absolutely continuous with respect to Lebesgue
measure. This can also be seen by observing that, if $L^a_t$ is the
local time of Brownian motion $B_t$ at $a$, then
\[
\int_0^t\varphi\circ u(B_s)u'^{\,2}(B_s)\,ds
 =\int_\RR\varphi\circ u(a)u'^{\,2}(a)L^a_t\,da,
\]
and applying the property stated in the introduction to $u$.

Nevertheless, the process $Y_t=\tilde u(B_t)$ is not in general a
semimartingale, since $\lsup{u}A_t$ is not of bounded variation when
$u$ is not a difference of convex functions (cf.~[5]). Thus the
occupation-time density property is extended to cases not covered by
[9].

\rem{Remark 10}For a general Dirichlet space under the assumptions of
Part II, the set of $u\in\DD$ satisfying the occupation-time density
property along the paths is stable under composition with
\pgo{163}%
one-variable Lipschitz functions. It clearly contains the functions
$u\in\DD$ that are semimartingales along the paths, a class which
contains differences of $p$-excessive functions belonging to $L^2(m)$
(cf.~[5]).

For Brownian motion with values in $\RR^d$, let $f\in H^1(\RR^d)$ have
canonical decomposition
\[
\tilde f(B_t)=\tilde f(B_0)
 +\int_0^t\bigl(\operatorname{grad}f(B_s),dB_s\bigr)
 +\lsup{f}A_t
\]
under $P_x$ for $x$ outside a polar set.

Let $\rho$ be a compactly supported measure that does not charge polar
sets. According to [4], for every $x$ and for $P_x$-almost every
$\omega$, the measure
\[
\zeta=\int_0^t\varepsilon_{B_s(\omega)}*\rho\,ds
\]
is absolutely continuous with respect to Lebesgue measure on $\RR^d$.
It follows from the property stated in the introduction that the
image under $f$ of the measure
$\operatorname{grad}^2f\,.\,\zeta$ is absolutely continuous with
respect to Lebesgue measure. We therefore obtain the following result.

\begin{enonce}{PROPOSITION 11}
Let $g\in L^2_{\mathrm{loc}}(\RR^d)$ be of the form
\[
g=f*\rho,
\]
where $f\in H^1(\RR^d)$ and $\rho$ is a compactly supported measure
that does not charge polar sets, so that $g\in H^1(\RR^d)$. Then $g$
satisfies quasi-everywhere the occupation-time density property along
the paths of $d$-dimensional Brownian motion.
\end{enonce}

\medskip
B. We now leave the Markovian setting and consider a probability space
$(\Omega,\Ftrib_t,\Ftrib,P)$ satisfying the usual conditions. We shall
study the occupation-time density property for Dirichlet processes,
i.e. processes of the form
\[
Y_t=Y_0+M_t+A_t,
\]
\pgo{164}%
where $M_t$ is a local martingale starting from zero and $A_t$ is a
process starting from zero with zero quadratic variation, in a sense to
be specified. In other words, we ask whether the image under
$s\mapsto Y_s(\omega)$ of the measure $d\cro{M^c,M^c}_s$ on $[0,t]$ is
absolutely continuous for $P$-almost every $\omega$.

Although several definitions are possible for processes with zero
quadratic variation, it is worth noting that the question is
fundamentally unaffected by an absolutely continuous change of
probability, by stopping, or by a time change.

Thus, from the property proved for one-dimensional Brownian motion in
Remark 9, one obtains the occupation-time density property for
processes of the form $u(X_t)$, where $u\in H^1(\RR)$ and $X_t$ is a
semimartingale that can be reduced to stopped Brownian motion by a time
change and an absolutely continuous change of probability.

This observation justifies our considering a continuous semimartingale
with decomposition
\begin{equation*}\tag{32}
X_t=X_0+N_t+B_t,
\end{equation*}
which does not need to satisfy $|dB_s|\ll d\cro{N,N}_s$, but which we take in
the space $H^2$ of semimartingales on $[0,1]$:
\begin{equation*}\tag{33}
\|X\|_{H^2}=\Bigl\|\,|X_0|+\cro{N,N}_1^{\frac12}
 +\int_0^1|dB_s|\,\Bigr\|_{L^2(\Omega,\Ftrib,P)}<\infty.
\end{equation*}

We adopt the following definitions.

\begin{enonce}{DEFINITION 12}
A process $Y_t$, $t\in[0,1]$, will be called a \uline{Dirichlet
process} if it can be written
\[
Y_t=Y_0+M_t+A_t,
\]
where $M_t$ is a martingale starting from zero such that
$\EE M_1^2<\infty$, and $A_t$ is a process starting from zero such that
\[
\EE\sum_{k=0}^{2^n-1}
 \bigl(A_{\frac{k+1}{2^n}}-A_{\frac{k}{2^n}}\bigr)^2
 \longrightarrow0 \qquad n\uparrow\infty.
\]
\end{enonce}
\pgo{165}%
The decomposition of $Y$ is then unique.

\begin{enonce}{DEFINITION 13}
We shall say that $Y$ satisfies the O.T.D. property (occupation-time density) if, $P$-a.s., the
image of $d\cro{M^c,M^c}_s$ on $[0,1]$ under
$s\mapsto Y_s$ is absolutely continuous with respect to Lebesgue
measure.
\end{enonce}

With the semimartingale $X$ satisfying (32)--(33), we associate the
seminorm $N^X$ defined by
\[
\bigl[N^X(f)\bigr]^2=\limsup_{n\uparrow\infty}
 \EE\sum_{k=0}^{2^n-1}
 \bigl[f(X_{\frac{k+1}{2^n}})-f(X_{\frac{k}{2^n}})\bigr]^2.
\]

\uline{We shall call the operating space associated with the
semimartingale} $X$ the space $\DD(X)$ of functions $f\in L^2(\RR)$
for which there exist infinitely differentiable compactly supported
functions $f_n$ ($f_n\in\Dsc$) such that
\[
N^X(f-f_n)+\|f-f_n\|_{L^2}\longrightarrow0
 \qquad n\uparrow\infty.
\]

This terminology is justified by the following proposition.

\begin{enonce}{PROPOSITION 14}
a) For every Borel function $f\in\DD(X)$, the process $f(X_t)$ is a
Dirichlet process whose martingale part is
$\int_0^t f^*(X_s)\,dN_s$ for a function $f^*$ satisfying
\[
\bigl[N^X(f)\bigr]^2=\int_\RR f^{*2}(a)\,\EE L^a_1\,da,
\]
where $L^a_t$ is the local time of the semimartingale $X$ at $a$. The
Dirichlet process $f(X_t)$ satisfies the O.T.D. property.

b) $\DD(X)$ contains every function $f\in L^2(\RR)$ of class $C^1$
whose derivative tends to zero at infinity.

c) If $\displaystyle\int_\alpha^\beta\frac{da}{\EE L^a_1}<+\infty$,
every function $f\in\DD(X)$ is equal almost everywhere on
$]\alpha,\beta[$ to an absolutely continuous function whose derivative
is equal to $f^*$ almost everywhere on $]\alpha,\beta[$.

d) If $\EE L^a_1\ge\lambda>0$ on $]\alpha,\beta[$, then the
restriction to $]\alpha,\beta[$ of every $f\in\DD(X)$ belongs to
$H^1(]\alpha,\beta[)$.
\end{enonce}
\pgo{166}%
\dem
1) Let $f$ be a Borel function in $\DD(X)$, and let $f_n\in\Dsc$ be
such that
\[
N^X(f-f_n)+\|f-f_n\|_{L^2}\to0.
\]

Since the semimartingale $(f_n-f_m)(X)$ belongs to the space $H^2$ of
semimartingales, we have (cf.~[12])
\begin{align*}
[N^X(f_n-f_m)]^2
 &=\limsup_{p\uparrow\infty}\EE\sum_{k=0}^{2^p-1}
   \bigl[(f_n-f_m)(X_{\frac{k+1}{2^p}})
   -(f_n-f_m)(X_{\frac{k}{2^p}})\bigr]^2\\
 &=\EE\int_0^1(f'_n-f'_m)^2(X_s)\,d\cro{N,N}_s\\
 &=\int_\RR(f'_n-f'_m)^2(a)\,\EE L^a_1\,da.
\end{align*}

Let $f^*$ be a Borel version of the limit of $f'_n$ in
$L^2\bigl((\EE L^a_1)\,da\bigr)$. If we write
\[
f_n(X_t)=f_n(X_0)+\int_0^t f'_n(X_s)\,dN_s+A^n_t,
\]
the stochastic integrals $\int_0^t f'_n(X_s)\,dN_s$ converge to
$\int_0^t f^*(X_s)\,dN$ in $L^2(\Omega,\Ftrib,P)$. Define $A_t$ by
\[
f(X_t)=f(X_0)+\int_0^t f^*(X_s)\,dN_s+A_t.
\]
Writing
\[
V_p(Z)=\sum_{k=0}^{2^p-1}
 \bigl(Z_{\frac{k+1}{2^p}}-Z_{\frac{k}{2^p}}\bigr)^2
 \qquad\text{for a process }Z,
\]
we have
\[
V_p(A)\le3V_p\bigl[(f-f_n)(X)\bigr]
 +3V_p\Bigl[\int_0^{\,\cdot}(f'_n-f^*)(X_s)\,dN_s\Bigr]
 +3V_p(B^n).
\]
Therefore
\[
\limsup_{p\uparrow\infty}\EE[\,V_p(A)\,]
 <3\bigl[N^X(f-f_n)\bigr]^2
 +3\|f'_n-f^*\|^2_{L^2(\EE L^a_1\,da)}.
\]
\pgo{167}%
The right-hand side can be made arbitrarily small by taking $n$ large
enough. Hence $f(X_t)$ is a Dirichlet process with martingale part
$\int_0^t f^*(X_s)\,dN_s$.

It also follows that
\[
\bigl[N^X(f)\bigr]^2=\EE\int_0^1 f^{*2}(X_s)\,dN_s
 =\int_\RR f^{*2}(a)\,\EE L^a_1\,da,
\]
and $f^*$ is unique up to equality $(\EE L^a_1)\,da$-almost
everywhere.

2) To prove the O.T.D. property for $f(X)$, we proceed as follows.

Consider the symmetric bilinear form $\Phi_\omega$ defined by
\[
\Phi_\omega(u,v)=\int_\RR u'(a)v'(a)L^a_1(\omega)\,da
 \qquad u,v\in\Dsc,
\]
and, setting $\nu=(\EE L^a_1)\,da$, the form $\Phi_\nu$ defined by
\[
\Phi_\nu(u,v)=\int_\RR u'(a)v'(a)\,\nu(da)
 \qquad u,v\in\Dsc.
\]

Since $a\mapsto L^a_1$ is càdlàg $P$-a.s., the forms $\Phi_\omega$ are
closable in $L^2(\RR)$ (cf.~[8]). It follows that $\Phi_\nu$ is
closable in $L^2(\RR)$ (cf.~[8], p.~45). Denote the associated
Dirichlet forms by $\overline\Phi_\omega$ and
$\overline\Phi_\nu$; they are regular, local, and conservative.

Let $u\in\Dsc\overline\Phi_\nu$, and let $u_n\in\Dsc$ satisfy
\[
\overline\Phi_\nu(u-u_n,u-u_n)+\|u-u_n\|_{L^2(\RR)}
 \xrightarrow[n\uparrow\infty]{}0.
\]
Then $u'_n$ is a Cauchy sequence in $L^2(\nu)$ and hence converges to
some $u^*\in L^2(\nu)$, with
\[
\overline\Phi_\nu(u,u)=\int u^{*2}(a)\,d\nu(a).
\]

Under these conditions, Lemmas 15 and 16 below show that, for
$\PP$-almost every $\omega$, $u\in\Dsc\overline\Phi_\omega$, and the
local energy measure of $u$ in
$(\Dsc\overline\Phi_\omega,\overline\Phi_\omega)$ is
\[
h\longrightarrow\int h(a)u^{*2}(a)L^a_t(\omega)\,da.
\]
\pgo{168}%
Part II, Corollary 6, then implies that the measure
\[
\varphi\longrightarrow
 \int\varphi\circ u(a)u^{*2}(a)L^a_t(\omega)\,da
\]
is absolutely continuous with respect to Lebesgue measure. Part 1)
shows that $\DD(X)\subset\Dsc\overline\Phi_\nu$ and that
$[N^X(f)]^2$ agrees with $\overline\Phi_\nu(f)$ for $f\in\DD(X)$.
This proves the asserted O.T.D. property.

3) To prove b), let $f\in L^2(\RR)$ be of class $C^1$ with derivative
in $C_0$. We may approximate $f$ in $L^2$ by functions $f_n\in\Dsc$
in such a way that $f'_n$ converges uniformly to $f'$. Then
\[
V_p\bigl[(f-f_n)(X)\bigr]
 \le\|f'-f'_n\|_\infty^2V_p(X),
\]
so that
\[
N^X(f-f_n)\le\|f'-f'_n\|_\infty\|X\|_{H^2},
\]
which proves the result.

4) Under assumption c), let $x,y\in]\alpha,\beta[$. With the notation
of part 1),
\begin{align*}
\Bigl|\int_x^y f^*(z)\,dz-f_n(y)+f_n(x)\Bigr|^2
 &\le\int_x^y\frac{da}{\EE L^a_1}
   \int_x^y\bigl(f^*(a)-f'_n(a)\bigr)^2\EE L^a_1\,da\\
 &\le\int_\alpha^\beta\frac{da}{\EE L^a_1}
   \|f^*-f'_n\|^2_{L^2(\nu)}.
\end{align*}
The result follows readily.

5) Under assumption d),
\[
\lambda\int_\alpha^\beta u'^{\,2}(a)\,da
 \le\int_\RR u'^{\,2}(a)\EE L^a_1\,da
 \qquad\forall u\in\Dsc.
\]
Thus, upon restriction to $]\alpha,\beta[$,
\[
\DD(X)\subset\Dsc(\overline\Phi_\nu)\subset H^1(]\alpha,\beta[).
\]

It remains to prove two lemmas.
\pgo{169}%
\begin{enonce}{LEMMA 15}
Let $\overline\Phi_\omega$ and $\overline\Phi_\nu$ be the regular
Dirichlet forms defined by
\begin{align*}
\overline\Phi_\omega(v,w)&=\int v'(a)w'(a)L^a_1(\omega)\,da
 \qquad v,w\in\Dsc,\\
\overline\Phi_\nu(v,w)&=\int v'(a)w'(a)\EE L^a_1\,da
 \qquad v,w\in\Dsc.
\end{align*}
If $u\in\Dsc\overline\Phi_\nu$, then
$u\in\Dsc\overline\Phi_\omega$ for almost every $\omega$.
\end{enonce}

\dem
If $u\in\Dsc\overline\Phi_\nu$, there exist $u_n\in\Dsc$ such that
$u_n\to u$ in $L^2(\RR)$ and
\[
\Phi_\nu(u_n-u_m,u_n-u_m)\to0
 \qquad\text{as }m,n\uparrow\infty.
\]
Choose a subsequence $n_i$ such that
\[
\sum_{i=1}^\infty
 \sqrt{\Phi_\nu(u_{n_{i+1}}-u_{n_i},u_{n_{i+1}}-u_{n_i})}<+\infty.
\]
Then, for $P$-almost every $\omega$,
\[
\sum_{i=1}^\infty
 \sqrt{\Phi_\omega(u_{n_{i+1}}-u_{n_i},u_{n_{i+1}}-u_{n_i})}<+\infty,
\]
and hence
\[
\Phi_\omega(u_{n_i}-u_{n_j},u_{n_i}-u_{n_j})\to0
 \qquad\text{as }n_i,n_j\uparrow\infty.
\]
Therefore $u\in\Dsc\overline\Phi_\omega$.

\begin{enonce}{LEMMA 16}
Under the same assumptions, let $u\in\Dsc\overline\Phi_\nu$, and let
$u_n\in\Dsc$ converge to $u$ in $L^2$ and with respect to
$\overline\Phi_\nu$.

Let $u^*$ be the limit of $u'_n$ in $L^2(\nu)$. Then, for almost every
$\omega$, the local energy measure of $u$ in
$\Dsc\overline\Phi_\omega$ is
\[
h\longrightarrow\int h(a)u^{*2}(a)L^a_1(\omega)\,da.
\]
\end{enonce}

\dem
The argument in the preceding proof shows that, for $\PP$-almost every
$\omega$, $u^*\in L^2(L^a_1(\omega)\,da)$ and
$u'_{n_i}\to u^*$ in $L^2(L^a_1(\omega)\,da)$.

Since the local energy measure $\mu_{u_{n_i}}$ of $u_{n_i}$ in
$(\Dsc\overline\Phi_\omega,\overline\Phi_\omega)$ is given by
\[
\langle\mu_{u_{n_i}},h\rangle
 =\int h(a)u_{n_i}'^{\,2}(a)L^a_1(\omega)\,da,
\]
the asserted result follows from the fact that, for every bounded
positive Borel function $h$, the mapping
$u\mapsto\langle\mu_u,h\rangle$ is continuous on
$(\Dsc\overline\Phi_\omega,
\overline\Phi_\omega+\|\cdot\|_{L^2})$ (cf.~[11]).
\pgo{170}%
\rem{Remark 17}\uline{If $f\in\DD(X)$ and $f$ is absolutely continuous
with derivative $f'$ in distributions sense, then one always has
$f^*=f'$ $(\EE L^a_1)\,da$-a.e.}

(Thus, in this case, the O.T.D. property for $f(X)$ follows immediately
from the property stated in the introduction and the O.T.D. property
for $X$.)

Indeed, let $f_n\in\Dsc(\RR)$ satisfy $N^X(f-f_n)\to0$ and
$\|f-f_n\|_{L^2}\to0$. Then $f'_n\to f^*$ in
$L^2(\EE L^a_1\,da)$ and $f'_n\to f'$ in the sense of distributions.
There is therefore a subsequence $n_i$ such that, for $\PP$-almost
every $\omega$, $f'_{n_i}\to f^*$ in
$L^2(L^a_1(\omega)\,da)$. But since $a\mapsto L^a_1(\omega)$ is
càdlàg, the set $\{L^a_1(\omega)>0\}$ differs from its interior by at
most a countable set. On that interior one must have $f^*=f'$
Lebesgue-a.e. Hence $f^*=f'$ $L^a_1(\omega)\,da$-a.e., and therefore
$(\EE L^a_1)\,da$-a.e.

\rem{Remark 18}\uline{Let $f\in H^1(\RR)$ be such that $f(X)$ is a
Dirichlet process. Then the martingale part of this Dirichlet process
is $\int_0^t f'(X_s)\,dN_s$ if and only if $f\in\DD(X)$.}

Indeed, if $f\in\DD(X)$, the result follows from the preceding remark.

Conversely, let $f\in H^1(\RR)$ be such that $f(X)$ is a Dirichlet
process with martingale part $\int_0^t f'(X_s)\,dN_s$.

First note that this martingale is well defined, since
\[
(X_1-a)^+=(X_0-a)^+
 +\int_0^1\mathbf{1}_{\{X_{s-}>a\}}\,dX_s+\tfrac12L^a_1
\]
implies
\[
\sup_a\EE L^a_1\le2\EE|X_1-X_0|+2\EE\int_0^1|dB_s|
 \le C\|X\|,
\]
and therefore
\[
\EE\Bigl(\int_0^t f'(X_s)\,dN_s\Bigr)^2
 =\int_\RR f'^{\,2}(a)\EE L^a_t\,da
 \le C\|X\|\,\|f'\|_{L^2}^2.
\]

Moreover, if $f_n\in\Dsc(\RR)$ satisfy
$\|f-f_n\|_{L^2}+\|f'-f'_n\|_{L^2}\to0$, then
\[
\bigl[N^X(f-f_n)\bigr]^2
 =\int_\RR(f'-f'_n)^2(a)\EE L^a_1\,da
 \le C\|X\|\,\|f'-f'_n\|_{L^2}^2,
\]
which implies that $f\in\DD(X)$.
\findem

\pgo{171}%
\newpage
\begin{center}
\uline{B I B L I O G R A P H Y}
\end{center}

\bigskip
\begin{list}{}{\setlength{\leftmargin}{3em}
 \setlength{\itemsep}{0.9em}\setlength{\labelwidth}{2.4em}}

\item[{[1]}] \uline{ANCONA A.}\\
-- Continuité des contractions dans les espaces de Dirichlet.
{\it Sém. th. du potentiel} No.~2, {\it Lect. Notes in Math.} No.~563,
Springer (1976).

\item[{[2]}] \uline{BLUMENTHAL R.M., GETOOR R.K.}\\
-- {\it Markov Processes and Potential Theory}. Academic Press (1968).

\item[{[3]}] \uline{BOULEAU N.}\\
-- Propriétés d'invariance du générateur étendu d'un processus de
Markov. Sém. Prob. XV, {\it Lect. Notes in Math.} 850, Springer (1981).

\item[{[4]}] \uline{BOULEAU N.}\\
-- Semi-martingales à valeurs $\RR^d$ et fonctions convexes.
{\it C.R. Acad. Sciences} Paris, vol.~292, pp.~87--90 (1981).

\item[{[5]}] \uline{CINLAR E., JACOD J., PROTTER P., SHARPE M.J.}\\
-- Semimartingales and Markov processes,
{\it Zeitschrift für Wahrscheinlichkeitstheorie} 54, 161--219 (1980).

\item[{[6]}] \uline{DELLACHERIE C., MEYER P.A.}\\
-- {\it Probabilités et potentiel, théorie des martingales},
Hermann, Paris (1980).

\item[{[7]}] \uline{DENY J., LIONS J.L.}\\
-- Les espaces de type Beppo Levi.
{\it Ann. Inst. Fourier} 5, 305--370 (1953/54).

\item[{[8]}] \uline{FUKUSHIMA M.}\\
-- {\it Dirichlet Forms and Markov Processes}. North-Holland (1980).

\item[{[9]}] \uline{GEMAN D., HOROWITZ J.}\\
-- Occupation densities. {\it Ann. of Probability}, Vol.~8, No.~1,
1--67 (1980).

\item[{[10]}] \uline{GETOOR R.K.}\\
-- Markov processes, Ray processes and right processes,
{\it Lect. Notes in Math.} 440, Springer (1975).

\item[{[11]}] \uline{LE JAN Y.}\\
-- Mesures associées à une forme de Dirichlet, applications.
{\it Bull. Soc. Math. France} 106, 61--112 (1978).

\item[{[12]}] \uline{LEPINGLE D.}\pgo{172}\\
-- La variation d'ordre $p$ des semi-martingales.
{\it Zeitschrift für Wahrscheinlichkeitstheorie} 36, 295--316 (1976).

\item[{[13]}] \uline{MEYER P.A.}\\
-- Un cours sur les intégrales stochastiques. Sém. Prob. X,
{\it Lect. Notes in Math.} 511, Springer (1976).

\item[{[14]}] \uline{MEYER P.A.}\\
-- La formule d'Itô pour le mouvement brownien d'après Brosamler.
Sém. Prob. XII, {\it Lect. Notes in Math.} 649, Springer (1978).

\item[{[15]}] \uline{REVUZ D.}\\
-- Mesures associées aux fonctionnelles additives de Markov I.
{\it Trans. Amer. Math. Soc.} 148, 501--531 (1970).

\end{list}

\vspace{3em}
\begin{center}
- o -
\end{center}

\vspace{3em}
\hspace*{0.5\textwidth}%
\begin{tabular}{@{}l@{}}
Nicolas BOULEAU\\
E.N.P.C. - C.E.R.M.A.\\
28, rue des Saints-Pères\\
75000 -- PARIS
\end{tabular}

\end{document}